\documentclass[12pt,reqno]{amsart}
\usepackage{amsmath,amssymb}
\usepackage[dvips]{graphicx,color}
\usepackage{latexsym,amssymb}
\usepackage{mathrsfs}
\usepackage[abbrev]{amsrefs}

\newtheorem{thm}{Theorem}[section]
\newtheorem{cor}[thm]{Corollaly}
\newtheorem{prop}[thm]{Proposition}
\newtheorem{lem}[thm]{Lemma}

\newtheorem{definition}[thm]{Definition}

\numberwithin{equation}{section}
\numberwithin{thm}{section}

\begin{document}

\begin{center}
 \textbf{\Large{Remark on the uniqueness of the mild solution \\ of SQG equation}}
\end{center}

\vskip5mm 
\centerline{Tsukasa Iwabuchi$^*$ and Ryoma Ueda$^{**}$}
\centerline{Mathematical Institute, Tohoku University}
\centerline{Sendai 980-8578 Japan}
\footnote[0]{\it{Mathematics Subject Classification}: 35Q35; 35Q86}
\footnote[0]{\it{Keywords}: quasi-geostrophic equation, mild solution, uniqueness} 
\footnote[0]{E-mail: $^*$t-iwabuchi@tohoku.ac.jp, $^{**}$ryoma.ueda.p6@dc.tohoku.ac.jp}

\begin{center}
 \begin{minipage}{120mm}
\small \textbf{Abstract.} We study the two-dimensional surface quasi-geostrophic equation. Motivated by the uniqueness for the three-dimensional incompressible Navier-Stokes equations, we demonstrate that the uniqueness of the mild solution of the two-dimensional surface quasi-geostrophic equation holds in the scaling critical Lebesgue space with a unique structure of the non-linear term.
 \end{minipage}
\end{center}

\section{Introduction.}
In this paper, we consider the surface quasi-geostrophic equation in $\mathbb R^2$.

\begin{equation}\tag{1}\label{eq1}
 \left\{
   \begin{aligned}
     &\partial_t \theta + u \cdot \nabla \theta - \Delta \theta = \ 0
           \quad &&\text{in} \  (0,T) \times \mathbb R^2, \\
     &u = \nabla^{\perp} \Lambda^{-1} \theta       
           \quad &&\text{in} \  (0,T) \times \mathbb R^2, \\
     &\theta(0,\cdot) = \ \theta_0(\cdot)
           \quad &&\text{in} \ \mathbb R^2,
   \end{aligned}
  \right.
\end{equation}
where $\nabla^{\perp} = (- \partial_{x_2}, \partial_{x_1}), \Lambda = (-\Delta)^{\frac{1}{2}}$. The equations are derived from general quasi-geostrophic equations in the special case of constant potential vorticity and buoyancy frequency (see~\cite{Pedlosky}). It is known that they describe the atmospheric motion and are useful in weather forecasting.  $\theta$ represents potential temperature, and $u$ denotes the atmospheric velocity vector. The purpose of this paper is to show the uniqueness of mild solutions of \eqref{eq1}.

\begin{definition}
\rm{Let} $T>0, \theta_0 \in L^2$. If $\theta : [0,T] \times \mathbb R^2 \to \mathbb R$ satisfies

\begin{equation}
 \left\{
   \begin{aligned}\notag
     & \theta \in C([0,T];L^2), \\
 & \theta(t) = e^{t\Delta}\theta_0 - \int_{0}^{t} \nabla \cdot e^{(t-s)\Delta} (u \theta) ds \quad \text{in} \ L^2 \ \text{for all} \quad t \in [0,T],
   \end{aligned}
  \right.
\end{equation}
then we call $\theta$ a mild solution of \eqref{eq1}.
\end{definition}

\begin{thm}\label{alpha=2}
 \rm{Let} $T>0$ and $\theta, \tilde{\theta} \in C([0,T];L^2)$ be mild solutions of \eqref{eq1} such that $\theta(0) = \tilde{\theta}(0)$. Then $\theta(t) = \tilde{\theta}(t)$ for all $t \in [0,T]$.
\end{thm}

Let us recall several known results about the uniqueness, where we consider $(-\Delta)^{\alpha/2}$ with $1 < \alpha \le 2$ instead of $-\Delta$ in the first equation of \eqref{eq1}. It is proved that the strong solution $\theta \in L^\infty(0,T;L^2) \cap L^2(0,T;H^{\frac{\alpha}{2}}) \cap L^q(0,T;L^p)$ is unique if $ \theta_0 \in L^2, p \ge 1, q >1, 1/p + \alpha/2q = \alpha/2 - 1/2 $ by Constantin and Wu \cite{MR1709781} . In regards to the uniqueness of mild solutions, Ferreira \cite{MR2836836} showed that the solution $\theta \in C([0,T];L^{\frac{2}{\alpha-1}})$ is unique if $1 < \alpha < 2$. The scale critical spaces of the two-dimensional surface quasi-geostrophic equation are same as those of the two-dimensional incompressible Navier-Stokes equations. It is a classical result that the weak solution $u \in L^{\infty}(0,T ; L^2) \cap L^2(0,T; \dot{H}^1)$ is unique if $u_0 \in L^2$ for the two-dimensional incompressible Navier-Stokes equations. We also refer to the paper \cite{MR4610908} by Cheskidov-Luo, who proved the non-uniqueness in $C([0, T] ; L^p)$ with $p<2$ for the two-dimensional incompressible Navier-Stokes equations. Uniqueness of the mild solution $u \in C([0,T];L^2)$ does not seem to be known. Regarding the three-dimensional incompressible Navier-Stokes equations, Meyer \cite{MR1724946}, Furioli, Lemari\'{e}-Rieusset, Terraneo \cite{MR1490408} and Monniaux \cite{MR1680809} demonstrated that the uniqueness of the mild solution $u \in C ([0, T ]; L^3 )$ holds. Lions-Masmoudi \cite{MR1876415} proved the uniqueness in $L^3$ by solving a dual problem. Furioli, Lemari\'{e}-Rieusset, Terraneo \cite{MR1490408} use Besov space and the proof of the main theorem in this paper is inspired by their result.

Furthermore, we have the following corollary since there is a global solution $\theta \in L^{\infty}(0,\infty ; L^2) \cap L^2(0,\infty; \dot{H}^1)$ for all $\theta_0 \in L^2$.

\begin{cor}
 \rm{Consider} the equations \eqref{eq1} in $(0, \infty) \times \mathbb R^2$. Then, for all $\theta_0 \in L^2$, there is a unique global mild solution $\theta \in C( [0, \infty) ; L^2)$.
\end{cor}

\vspace{5mm}
\textbf{Notation.} Let $\{ \phi_j \}_{j \in \mathbb Z}$ be the Littlewood-Paley dyadic decomposition such that $ \phi_j $ is a non-negative function in $C_0^\infty$ and satisfies that
\begin{align*}
 &\text{supp} \ \phi_0 \subset \big\{ \xi \in \mathbb R^2 \ \big| \ 2^{-1} \le \lvert \xi \rvert \le 2 \big\}, \quad \phi_j(\xi) = \phi_0 \Bigl( \frac{\xi}{2^j} \Bigr) \ \text{for any} \ \xi \in \mathbb R^2,\\
 &\sum_{j \in \mathbb Z} \phi_j(\xi) = 1 \ \text{for any} \ \xi \in \mathbb R^2 \setminus \{0\}.
\end{align*}
We define homogeneous Besov spaces ${\dot{B}}_{p, q}^{s} = {\dot{B}}_{p, q}^{s}(\mathbb R^2)$ for $s \in \mathbb R, 1 \le p, q \le \infty$ as follows.
\begin{align*}
 {\dot{B}}_{p, q}^{s} = {\dot{B}}_{p, q}^{s}(\mathbb R^2) := \big\{f \in \mathcal S' / \mathcal P \ \big| \ \bigl\Vert f \bigr\Vert_{{\dot{B}}_{p, q}^{s}(\mathbb R^2)}< \infty \big\},
\end{align*}
where $\mathcal S'$ is the dual space of Schwartz class $\mathcal S$, $\mathcal P$ is the set of all polynominals of 2 real variables and
\begin{align*}
 \bigl\Vert f \bigr\Vert_{{\dot{B}}_{p, q}^{s}(\mathbb R^2)} := \Bigl\Vert \Big\{ 2^{sj} \bigl\Vert \phi_j(D) f \bigr\Vert_{L^p(\mathbb R^2)} \Big\}_{j \in \mathbb Z} \Bigr\Vert_{l^q(\mathbb Z)}.
\end{align*}
We denote $f_j := \phi_j(D)f = {\mathcal F}^{-1}[\phi_j(\xi) \hat{f}(\xi)]$ for $j \in \mathbb Z$ and introduce a notation for nonlinear term in the integral equation, 
\begin{align*}
 \theta(t) = e^{t\Delta}\theta_0 - \int_{0}^{t} \nabla \cdot e^{(t-s)\Delta} (u \theta) ds =: e^{t\Delta}\theta_0 + B(u\theta)(t).
\end{align*}

\section{Preliminaries.}
In this section, we introduce some lemmas for the proof. The first four lemmas are elemental and we refer to books by Triebel \cite{MR0781540}, Danchin \cite{MR2768550} and Grafakos \cite{MR3243741}.

\begin{lem}\label{LP estimate} 
\rm{Let} $s \in \mathbb R, \alpha >0, 1 \le p, q \le \infty$. Then
\begin{align}
  &f = \sum_{j \in \mathbb Z} \phi_j(D)f \ \text{in} \ \mathcal S' \ \text{for all} \ f \in L^2 \notag, \\
  &\sup_{j \in \mathbb Z} \lVert \phi_j(D) \rVert_{L^p \rightarrow L^p} < \infty ,\notag \\
  &\bigl\Vert \phi_j(D) \nabla f \bigr\Vert_{L^p} \le C 2^j  \bigl\Vert \phi_j(D) f \bigr\Vert_{L^p} \ \text{for all} \ f \in L^p, \ \text{all} \ j \in \mathbb Z \notag, \\
  &\bigl\Vert e^{t\Delta} f \bigr\Vert_{{\dot{B}}_{p, q}^{s}} \le C t^{-\frac{\alpha}{2}} \bigl\Vert f \bigr\Vert_{{\dot{B}}_{p, q}^{s-\alpha}} \ \text{for all} \ f \in {\dot{B}}_{p, q}^{s-\alpha} \label{smoothing effect}, \\
  &\bigl\Vert \phi_j(D) e^{-t\Delta} f \bigr\Vert_{L^p} \le C e^{-Ct{2^{2j}}} \bigl\Vert \phi_j(D) f \bigr\Vert_{L^p} \ \text{for all} \ f \in L^p, \ \text{all} \ j \in \mathbb Z \notag, \\
  &\int_{0}^{\infty} \bigl\Vert e^{-t\Delta} f \bigr\Vert_{{\dot{B}}_{\infty, 1}^{s}} dt \le C \bigl\Vert f \bigr\Vert_{{\dot{B}}_{\infty, 1}^{s-2}} \ \text{for all} \ f \in {\dot{B}}_{\infty, 1}^{s-2}\label{maximum regularity}.
\end{align}
\end{lem}

\begin{lem}\label{Besov dual}
\rm{Let} $1 \le p, q < \infty, s \in \mathbb R$ and $p', q'$ satisfy $1 = 1/p + 1/p' = 1/q +1/q'$. Then
\begin{equation}\notag
 \begin{aligned}
  ({\dot{B}}_{p,q}^{s})^{*} = {\dot{B}}_{p',q'}^{-s}, \quad
  \bigl\Vert f \bigr\Vert_{{\dot{B}}_{p,q}^{s}} = \sup_{\Vert \varphi \Vert_{{\dot{B}}_{p',q'}^{-s}} = 1} \bigl\vert \langle f,\varphi \rangle \bigr\vert \ \text{for} \ f \in {\dot{B}}_{p,q}^{s} .
 \end{aligned}
\end{equation}
\end{lem}

\begin{lem}\label{Besov-Holder}
\rm{Let} $s \in \mathbb R, 1 \le p, p_1, p_2, p_3, p_4, q \le \infty, 1/p = 1/p_1 + 1/p_2 = 1/p_3 + 1/p_4$. Then
\begin{align*}\notag
 &\bigl\Vert fg \bigr\Vert_{{\dot{B}}_{p,q}^{s}} \le C\Bigl( \bigl\Vert f \bigr\Vert_{{\dot{B}}_{p_1,q}^{s}} \bigl\Vert g \bigr\Vert_{L^{p_2}} + \bigl\Vert f \bigr\Vert_{L^{p_3}}\bigl\Vert g \bigr\Vert_{{\dot{B}}_{p_4,q}^{s}}\Bigr) \ \text{for} \ f \in {\dot{B}}_{p_1,q}^{s} \cap L^{p_3}, g \in L^{p_2} \cap {\dot{B}}_{p_4,q}^{s}, \\
 &\biggl\Vert \sum_{k \ge l+3}f_k g_l  \biggr\Vert_{{\dot{B}}_{p,q}^{s}} \le 
  C \bigl\Vert f \bigr\Vert_{{\dot{B}}_{p_1,q}^{0}} \bigl\Vert g \bigr\Vert_{{\dot{B}}_{p_2,q}^{s}} \ \text{for} \ f \in {\dot{B}}_{p_1,q}^{0}, g \in {\dot{B}}_{p_2,q}^{s}.
\end{align*}
\end{lem}

\begin{lem}\label{Besov embedding}
\rm{Let} $s_1, s_2 \in \mathbb R, s_2 < s_1, 1 \le p_1, p_2, q \le \infty, p_1 < p_2, s_1 - 2/p_1 = s_2 - 2/p_2$. Then ${\dot{B}}_{p_1,q}^{s_1} \hookrightarrow {\dot{B}}_{p_2,q}^{s_2}$, i.e.,
\begin{align}\notag
 \bigl\Vert f \bigr\Vert_{{\dot{B}}_{p_2, q}^{s_2}} \le C \bigl\Vert f \bigr\Vert_{{\dot{B}}_{p_1, q}^{s_1}} \ \text{for} \ f \in {\dot{B}}_{p_1,q}^{s_1}.
\end{align}
\end{lem}

We have the estimate of non-linear term by next lemma.

\begin{lem}\label{Fourier}
\rm{Let} $s>0, 1 \le p, p_1, p_2 < \infty, 1/p = 1/p_1 + 1/p_2, 1 \le q \le \infty$ and $\alpha$ be a multi-index with $\lvert \alpha \rvert \le 2$. Suppose that a bilinear Fourier multiplier $m(\cdot, \cdot)$ satisfies that
\begin{align*}
  \lvert \partial_{\xi-\eta, \eta}^{\alpha} m(\xi-\eta, \eta) \rvert \le \frac{C_\alpha}{\big(\lvert \xi-\eta \rvert + \lvert \eta \rvert \big)^{\lvert \alpha \rvert}}.
\end{align*}
Then
\begin{align*} 
  &\Biggl\Vert {\mathcal F}^{-1} \bigg[\int_{{\mathbb R}^2} m(\xi - \eta, \eta) \sum_{\lvert k-l \rvert \le 2} \hat{f_k}(\xi - \eta) \hat{g_l}(\eta)d\eta \bigg] \Biggr\Vert_{{\dot{B}}_{p,q}^{s}} \le C \bigl\Vert f \bigr\Vert_{{\dot{B}}_{p_1,q}^{s}} \bigl\Vert g \bigr\Vert_{L^{p_2}} \\
  &\text{for all} \ f \in {\dot{B}}_{p_1, q}^{s} \ \text{and} \ g \in L^{p_2}.\end{align*}
\end{lem}

\noindent
\textbf{Proof.} 
We write
\begin{align*}
 &\Biggl\Vert {\mathcal F}^{-1} \bigg[\int_{{\mathbb R}^2} m(\xi - \eta, \eta) \sum_{\lvert k-l \rvert \le 2} \hat{f_k}(\xi - \eta) \hat{g_l}(\eta)d\eta \bigg] \Biggr\Vert_{{\dot{B}}_{p,q}^{s}} \\
 = &\Bigg\{ \sum_{j \in \mathbb Z} \Bigg(2^{sj} \biggl\Vert \phi_j(D) {\mathcal F}^{-1} \bigg[\int_{{\mathbb R}^2} m(\xi - \eta, \eta) \sum_{\lvert k-l \rvert \le 2} \hat{f}_k(\xi - \eta) \hat{g}_l(\eta)d\eta \bigg] \biggr\Vert_{L^p} \Bigg)^q \Bigg\}^{\frac{1}{q}} \\
 \le &C \Bigg\{ \sum_{j \in \mathbb Z} \Bigg(2^{sj} \sum_{k \ge j-5} \sum_{\lvert l-k \rvert \le 2} \biggl\Vert {\mathcal F}^{-1} \bigg[\int_{{\mathbb R}^2} m(\xi - \eta, \eta) \hat{f}_k(\xi - \eta) \hat{g}_l(\eta)d\eta \bigg] \biggr\Vert_{L^p} \Bigg)^q \Bigg\}^{\frac{1}{q}}. \\
\end{align*}
By the boundedness of $m(D_1, D_2)$ in $L^p$ (see~\cite{MR3028565}), 
\begin{align*}
 &\Biggl\Vert {\mathcal F}^{-1} \bigg[\int_{{\mathbb R}^2} m(\xi - \eta, \eta) \sum_{\lvert k-l \rvert \le 2} \hat{f_k}(\xi - \eta) \hat{g_l}(\eta)d\eta \bigg] \Biggr\Vert_{{\dot{B}}_{p,q}^{s}}  \\
 \le & C \Bigg\{ \sum_{j \in \mathbb Z} \Bigg( 2^{sj} \sum_{k \ge j-5} \sum_{\lvert l-k \rvert \le 2} \bigl\Vert f_k \bigr\Vert_{L^{p_1}} \bigl\Vert g_l \bigr\Vert_{L^{p_2}} \Bigg)^q \Bigg\}^{\frac{1}{q}} \\
 \le & C \bigl\Vert f \bigr\Vert_{{\dot{B}}_{p_1,q}^{s}} \bigl\Vert g \bigr\Vert_{L^{p_2}},
\end{align*} 
which completes the proof. \hfill $\square$ \\

The next two propsitions are also basic. Proposition~\ref{linear} is well-known. Proposition~\ref{nonlinear} holds since $\theta$ is a mild solution of \eqref{eq1}.

\begin{prop}\label{linear}
\rm{If} $2 < p \le \infty$, \rm{then}
\begin{equation}\notag
 \lim_{t \rightarrow 0} t^{\frac{1}{2}-\frac{1}{p}} \big\lVert e^{t \Delta} \theta_0 \big\rVert_{L^p} = 0 \ \text{for all} \ \theta_0 \in L^2.
\end{equation}
\end{prop}

\begin{prop}\label{nonlinear}
\rm{Suppose that} $\theta$ is a mild solution of \eqref{eq1}. Then
\begin{equation}\notag
 \lim_{t \rightarrow 0} \big\lVert \theta(t) - e^{t\Delta} \theta_0 \big\rVert_{L^2}= 0.
\end{equation}
\end{prop}

The following two propositions are essential in considering the uniqueness of solutions in $L^2$ by the estimates related to the Besov norm.
\begin{prop}
 \rm{Suppose that} $1<p<2$ and $\theta$ is a mild solution of \eqref{eq1}. Then
\begin{equation}\notag
 \bigl\Vert \theta(t) - e^{t\Delta} \theta_0 \bigr\Vert_{{\dot{B}}_{p,\infty}^{-1+\frac{2}{p}}} \le C \big\lVert \theta \big\rVert_{L^\infty(0,T;L^2)}^2 \ \text{for} \ t \in [0,T].
\end{equation} 
\end{prop}

\noindent
\textbf{Proof.}
Since $\theta$ satisfies the integral equation of \eqref{eq1} and  $\Big({\dot{B}}_{p,\infty}^{-1+\frac{2}{p}}\Big)^{*} = {\dot{B}}_{p',1}^{1-\frac{2}{p}}$,
\begin{align*}
 \bigl\Vert \theta - e^{t\Delta} \theta_0 \bigr\Vert_{{\dot{B}}_{p,\infty}^{-1+\frac{2}{p}}} 
 &= \sup_{\Vert \varphi \Vert_{{\dot{B}}_{p',1}^{1-\frac{2}{p}}} = 1} \Biggl\vert \bigg\langle \int_{0}^{t} \nabla \cdot e^{(t-s)\Delta} (u \theta) ds,\varphi \bigg\rangle \Biggr\vert \\
 &= \sup_{\Vert \varphi \Vert_{{\dot{B}}_{p',1}^{1-\frac{2}{p}}} = 1} \Biggl\vert \int_{0}^{t} \Big\langle u \theta, e^{(t-s)\Delta} \nabla \varphi \Big\rangle ds \Biggr\vert \\
 &\le \sup_{\Vert \varphi \Vert_{{\dot{B}}_{p',1}^{1-\frac{2}{p}}} = 1}  \int_{0}^{t} \Big\lvert\Big\langle u \theta, e^{(t-s)\Delta} \nabla \varphi \Big\rangle \Big\rvert ds.
\end{align*}
By applying the \rm{H}$\ddot{\text{o}}$\rm{lder} inequality, boundedness of Riesz transformation and the embedding ${\dot{B}}_{\infty,1}^{0} \hookrightarrow L^\infty$,
\begin{align*}
 \bigl\Vert \theta - e^{t\Delta} \theta_0 \bigr\Vert_{{\dot{B}}_{p,\infty}^{-1+\frac{2}{p}}} 
  &\le \sup_{\Vert \varphi \Vert_{{\dot{B}}_{p',1}^{1-\frac{2}{p}}} = 1}  \int_{0}^{t} \bigl\Vert u \theta \bigr\Vert_{L^1} \bigl\Vert e^{(t-s)\Delta} \nabla \varphi \bigr\Vert_{L^\infty} ds \\
  &\le C \sup_{\Vert \varphi \Vert_{{\dot{B}}_{p',1}^{1-\frac{2}{p}}} = 1}  \int_{0}^{t} \big\lVert \theta \big\rVert_{L^2}^2 \bigl\Vert e^{(t-s)\Delta} \nabla \varphi \bigr\Vert_{{\dot{B}}_{\infty,1}^{0}} ds \\
  &\le C \big\lVert \theta \big\rVert_{L^\infty(0,T;L^2)}^2 \sup_{\Vert \varphi \Vert_{{\dot{B}}_{p',1}^{1-\frac{2}{p}}} = 1}  \int_{0}^{t}  \bigl\Vert e^{(t-s)\Delta} \nabla \varphi \bigr\Vert_{{\dot{B}}_{\infty,1}^{0}} ds.
\end{align*}
By maximum regularity estimate \eqref{maximum regularity} and the embedding ${\dot{B}}_{p',1}^{1-\frac{2}{p}} \hookrightarrow {\dot{B}}_{\infty,1}^{-1}$,
\begin{align*}
 \bigl\Vert \theta - e^{t\Delta} \theta_0 \bigr\Vert_{{\dot{B}}_{p,\infty}^{-1+\frac{2}{p}}}
  &\le C \big\lVert \theta \big\rVert_{L^\infty(0,T;L^2)}^2 \sup_{\Vert \varphi \Vert_{{\dot{B}}_{p',1}^{1-\frac{2}{p}}} = 1} \big\lVert \varphi \big\rVert_{{\dot{B}}_{\infty,1}^{-1}}& \\
  &\le C \big\lVert \theta \big\rVert_{L^\infty(0,T;L^2)}^2 \sup_{\Vert \varphi \Vert_{{\dot{B}}_{p',1}^{1-\frac{2}{p}}} = 1} \big\lVert \varphi \big\rVert_{{\dot{B}}_{p',1}^{1-\frac{2}{p}}}& \\
  &= C \big\lVert \theta \big\rVert_{L^\infty(0,T;L^2)}^2 ,
 \end{align*}
which completes the proof. \hfill $\square$ \\

\begin{prop}\label{shiage}
\rm{If} $ f,g \in L^2, f=g$ in ${\dot{B}}_{2,\infty}^{0}$, then $f=g$ in $L^2$.
\end{prop}

\noindent
\textbf{Proof.}
We can write
\begin{align*}
 &{\dot{B}}_{2,\infty}^{0} = \bigg\{ f \in \mathcal S' \ \bigg| \ f = \sum_{j \in \mathbb Z} \phi_j(D)f \ \text{in} \ \mathcal S', \Vert f \Vert_{{\dot{B}}_{2, \infty}^{0}} < \infty \bigg\}, \\
 &L^2 =  \bigg\{ f \in \mathcal S' \ \bigg| \ f = \sum_{j \in \mathbb Z} \phi_j(D)f \ \text{in} \ \mathcal S', \Vert f \Vert_{L^2} < \infty \bigg\}.
\end{align*}
Since $f=g$ in ${\dot{B}}_{2,\infty}^{0}$, $f=g$ in $\mathcal S'$. Given this result and the assumption that $f, g \in L^2$, it follows that $f=g \ a.e.$ Consequently, $f=g$ in $L^2$. \hfill $\square$ \\

\section{Proof of theorem.} 
\noindent
Step 1. Let $\psi := \theta - \tilde{\theta}, \theta_0 := \theta(0) = \tilde\theta(0), \tilde{u} := \nabla^\perp \Lambda^{-1} \tilde{\theta}, 1 < p < 2 $. We show that $\big\lVert \psi \big\rVert_{L^2} = 0$ in $[0,\delta]$ for some small $\delta>0$. To this end, we consider $\big\lVert \psi \big\rVert_{{\dot{B}}_{p,\infty}^{-1+\frac{2}{p}}}$ instead of the $L^2$ norm. We start by rewriting
\begin{equation}\notag
 \psi = - \int_{0}^{t} \nabla \cdot e^{(t-s)\Delta} (u \theta - \tilde{u} \tilde{\theta}) ds
\end{equation}
 with
\begin{equation}\notag
 \begin{split}
 u \theta - \tilde{u} \tilde{\theta}
 = \frac{1}{2} \{ (\nabla^{\perp} \Lambda^{-1} \psi) \theta + (\nabla^{\perp} \Lambda^{-1} \theta) \psi + (\nabla^{\perp} \Lambda^{-1} \psi) \tilde{\theta} + (\nabla^{\perp} \Lambda^{-1} \tilde{\theta}) \psi \}.
 \end{split}
\end{equation}
This implies that
\begin{align*}
 \bigl\Vert \psi \bigr\Vert_{{\dot{B}}_{p,\infty}^{-1+\frac{2}{p}}} 
 \le& \ \frac{1}{2} \Bigg\{ \Biggl\Vert {\int_{0}^{t} \nabla \cdot e^{(t-s)\Delta}  \big\{ (\nabla^{\perp} \Lambda^{-1} \psi) \theta + (\nabla^{\perp} \Lambda^{-1} \theta) \psi \big\} ds} \Biggr\Vert_{{\dot{B}}_{p,\infty}^{-1+\frac{2}{p}}} \\
 & \ + \Biggl\Vert {\int_{0}^{t} \nabla \cdot e^{(t-s)\Delta}  \big\{(\nabla^{\perp} \Lambda^{-1} \psi) \tilde{\theta} + (\nabla^{\perp} \Lambda^{-1} \tilde{\theta}) \psi \big\} ds} \Biggr\Vert_{{\dot{B}}_{p,\infty}^{-1+\frac{2}{p}}} \Bigg\} \\
 =&: \frac{1}{2}(I + J).
\end{align*}
We only consider \textit{I} from now on, since the estimate regarding \textit{I} can also be applied to \textit{J}.
Substituting $\theta(t) = e^{t\Delta} \theta_0 + B(u \theta)(t)$,
\begin{align*}
 I \le& \Biggl\Vert {\int_{0}^{t} \nabla \cdot e^{(t-s)\Delta}  \big\{ (\nabla^{\perp} \Lambda^{-1} \psi) e^{s\Delta}\theta_0 + (\nabla^{\perp} \Lambda^{-1} (e^{s\Delta}\theta_0)) \psi \big\} ds} \Biggr\Vert_{{\dot{B}}_{p,\infty}^{-1+\frac{2}{p}}} \\
 & \ + \Biggl\Vert {\int_{0}^{t} \nabla \cdot e^{(t-s)\Delta}  \big\{(\nabla^{\perp} \Lambda^{-1} \psi) B(u \theta)(s) + (\nabla^{\perp} \Lambda^{-1} B(u \theta)(s)) \psi \big\} ds} \Biggr\Vert_{{\dot{B}}_{p,\infty}^{-1+\frac{2}{p}}} \\
 =&: I_{\text{linear}} + I_{\text{nonlinear}}.
\end{align*}

\noindent 
By Bony's decomposition (see~\cite{MR0631751}), we write 
\begin{align*}
 I_{\text{linear}} 
 &\le  \Biggl\Vert \int_{0}^{t} \nabla \cdot e^{(t-s)\Delta}  \sum_{k \ge l+3} \Big\{(\nabla^{\perp} \Lambda^{-1} \psi_k) (e^{s\Delta}\theta_0)_l + (\nabla^{\perp} \Lambda^{-1} (e^{s\Delta}\theta_0)_l) \psi_k \Big\} ds \Biggr\Vert_{{\dot{B}}_{p,\infty}^{-1+\frac{2}{p}}} \\
 &\hspace{8mm}+ \Biggl\Vert \int_{0}^{t} \nabla \cdot e^{(t-s)\Delta}  \sum_{\lvert k-l \rvert \le 2} \Big\{(\nabla^{\perp} \Lambda^{-1} \psi_k) (e^{s\Delta}\theta_0)_l + (\nabla^{\perp} \Lambda^{-1} (e^{s\Delta}\theta_0)_l) \psi_k \Big\} ds \Biggr\Vert_{{\dot{B}}_{p,\infty}^{-1+\frac{2}{p}}} \\
 &\hspace{8mm}+ \Biggl\Vert \int_{0}^{t} \nabla \cdot e^{(t-s)\Delta} \sum_{l \ge k+3} \Big\{(\nabla^{\perp} \Lambda^{-1} \psi_k) (e^{s\Delta}\theta_0)_l + (\nabla^{\perp} \Lambda^{-1} (e^{s\Delta}\theta_0)_l) \psi_k \Big\} ds \Biggr\Vert_{{\dot{B}}_{p,\infty}^{-1+\frac{2}{p}}} \\
 &=: I_{\text{linear}}^1 + I_{\text{linear}}^2 + I_{\text{linear}}^3,
\end{align*}
\begin{align*}
 I_{\text{nonlinear}}
  &\le \Biggl\Vert {\int_{0}^{t} \nabla \cdot e^{(t-s)\Delta}  \sum_{k \ge l+3} \Big\{(\nabla^{\perp} \Lambda^{-1} \psi_k) (B(u \theta))_l + (\nabla^{\perp} \Lambda^{-1} (B(u \theta))_l) \psi_k \Big\} ds} \Biggr\Vert_{{\dot{B}}_{p,\infty}^{-1+\frac{2}{p}}} \\
  &\hspace{8mm} +\Biggl\Vert {\int_{0}^{t} \nabla \cdot e^{(t-s)\Delta}  \sum_{\lvert k-l \rvert \le 2} \Big\{(\nabla^{\perp} \Lambda^{-1} \psi_k) (B(u \theta))_l + (\nabla^{\perp} \Lambda^{-1} (B(u \theta))_l) \psi_k \Big\} ds} \Biggr\Vert_{{\dot{B}}_{p,\infty}^{-1+\frac{2}{p}}} \\
  & \hspace{8mm} + \Biggl\Vert {\int_{0}^{t} \nabla \cdot e^{(t-s)\Delta}  \sum_{l \ge k+3} \Big\{(\nabla^{\perp} \Lambda^{-1} \psi_k) (B(u \theta))_l + (\nabla^{\perp} \Lambda^{-1} (B(u \theta))_l) \psi_k \Big\} ds} \Biggr\Vert_{{\dot{B}}_{p,\infty}^{-1+\frac{2}{p}}} \\
  &=: I_{\text{nonlinear}}^1 + I_{\text{nonlinear}}^2 + I_{\text{nonlinear}}^3.
\end{align*}

First, we consider $I_{\text{linear}}^1, I_{\text{linear}}^2$ and $I_{\text{linear}}^3$. Especially, we only handle the estimate of $I_{\text{linear}}^2$, because it is easier to check the following inequality: 
\begin{align}\label{liner I_1,I_3}
I_{\text{linear}}^1, I_{\text{linear}}^3 \le  C \Big(\sup_{s \in [0,\delta]} s^{\frac{1}{2}-\frac{1}{2p}} \lVert e^{s\Delta}\theta_0 \rVert_{L^{2p}} \Big) \Big(\sup_{s \in [0,\delta]} \lVert \psi \rVert_{{\dot{B}}_{p,\infty}^{-1 + \frac{2}{p}}} \Big).
\end{align}
Actually, one can get these inequalities, applying the facts about smoothing effect of $e^{t\Delta}$ \eqref{smoothing effect}, embedding property of homogeneous Besov spaces (Lemma~\ref{Besov embedding}) and the product estimate of Besov norm (Lemma~\ref{Besov-Holder}). \\

Now, we write
\begin{align}\label{fg}
 (\nabla^{\perp} \Lambda^{-1} f) g + (\nabla^{\perp} \Lambda^{-1} g) f 
 &=  {\mathcal F}^{-1} \bigg[\int_{{\mathbb R}^2} \bigg(\frac{(\xi - \eta)^{\perp}}{\lvert \xi - \eta \rvert} + \frac{\eta^\perp}{\lvert \eta \rvert} \bigg) \hat{f}(\xi - \eta) \hat{g}(\eta)d\eta \bigg] \notag \\
 &= {\mathcal F}^{-1} \bigg[\int_{{\mathbb R}^2} \frac{\xi^\perp}{\lvert \xi - \eta \rvert} \hat{f}(\xi - \eta) \hat{g}(\eta)d\eta \bigg] \notag \\
   &\hspace{8mm}+ {\mathcal F}^{-1} \bigg[\int_{{\mathbb R}^2} \frac{\xi \cdot (\xi - 2\eta)}{\lvert \xi - \eta \rvert + \lvert \eta \rvert} \cdot \frac{1}{\lvert \xi - \eta \rvert}\hat{f}(\xi - \eta) \cdot \frac{\eta^\perp}{\lvert \eta \rvert}\hat{g}(\eta)d\eta \bigg] \notag \\
  &= \nabla^\perp \{(\Lambda^{-1} f) g \} + \nabla m(D_1, D_2)\{(\Lambda^{-1} f)(\nabla^\perp \Lambda^{-1} g)\},
\end{align}
where $m(D_1, D_2)(fg) := \displaystyle{ {\mathcal F}^{-1} \bigg[\int_{{\mathbb R}^2} \frac{\xi - 2\eta}{\lvert \xi - \eta \rvert + \lvert \eta \rvert} \hat{f}(\xi - \eta) \hat{g}(\eta)d\eta \bigg]}$. 
We replace $f$ and $g$ in \eqref{fg} with $\psi_k$ and $(e^{s\Delta} \theta_0)_l$ respectively. We estimate $I_{\text{linear}}^2$ as follows: 
\begin{align*}
 I_{\text{linear}}^2 &\le \Biggl\Vert \int_{0}^{t} \nabla \cdot e^{(t-s)\Delta}  \sum_{\lvert k-l \rvert \le 2} \nabla m(D_1, D_2) \{(\Lambda^{-1} \psi_k)(\nabla^\perp \Lambda^{-1} (e^{s\Delta}\theta_0)_l)\} ds  \Biggr\Vert_{{\dot{B}}_{p,\infty}^{-1+\frac{2}{p}}} \\
 &\le C \Biggl\Vert \int_{0}^{t} e^{(t-s)\Delta}  \sum_{\lvert k-l \rvert \le 2} m(D_1, D_2) \{(\Lambda^{-1} \psi_k)(\nabla^\perp \Lambda^{-1} (e^{s\Delta}\theta_0)_l)\} ds  \Biggr\Vert_{{\dot{B}}_{p,\infty}^{1+\frac{2}{p}}},
\end{align*} 
since $\nabla \cdot \nabla^\perp \{(\Lambda^{-1} f) g \} = 0$ for the first term of \eqref{fg}.
By smoothing effect \eqref{smoothing effect}, Lemma~\ref{Fourier}, the embedding and the boundedness of Riesz transformation, 
\begin{align}\label{linear I_2}
 I_{\text{linear}}^2 &\le C \int_{0}^{t} (t-s)^{-\frac{1}{2}(1+\frac{1}{p})} \biggl\Vert  \sum_{\lvert k-l \rvert \le 2} m(D_1, D_2) (\Lambda^{-1} \psi_k)(\nabla^\perp \Lambda^{-1} (e^{s\Delta}\theta_0)_l) \biggr\Vert_{{\dot{B}}_{p,\infty}^{\frac{1}{p}}} ds \notag \\
 &\le C \int_{0}^{t} (t-s)^{-\frac{1}{2}(1+\frac{1}{p})} \Bigl\Vert \Lambda^{-1} \psi \Big\rVert_{{\dot{B}}_{2p,\infty}^{\frac{1}{p}}} \Big\lVert \nabla^\perp \Lambda^{-1} e^{s\Delta}\theta_0 \Bigr\Vert_{L^{2p}} ds \notag \\
 &\le C \Big(\sup_{s \in [0,\delta]} \lVert \psi \rVert_{{\dot{B}}_{p,\infty}^{-1 + \frac{2}{p}}}\Big) \Big(\sup_{s \in [0,\delta]} s^{\frac{1}{2}-\frac{1}{2p}} \lVert e^{s\Delta}\theta_0 \rVert_{L^{2p}} \Big) .
\end{align} 

Next, we consider $I_{\text{nonlinear}}^1, I_{\text{nonlinear}}^2$ and $I_{\text{nonlinear}}^3$. By duality in Besov spaces, we have
\begin{align*}
 I_{\text{nonlinear}}^1  &\le \sup_{\Vert \varphi \Vert_{{\dot{B}}_{p',1}^{1-\frac{2}{p}}} = 1}  \int_{0}^{t} \Biggl\Vert \sum_{k \ge l+3} \Big\{(\nabla^{\perp} \Lambda^{-1} \psi_k) (B(u \theta))_l + (\nabla^{\perp} \Lambda^{-1} (B(u \theta))_l) \psi_k \Big\}\Biggr\Vert_{{\dot{B}}_{p,\infty}^{-2+\frac{2}{p}}} \\ 
 &\hspace{5cm} \times \big\lVert e^{(t-s)\Delta} \nabla \varphi \big\rVert_{{\dot{B}}_{p',1}^{2-\frac{2}{p}}} ds,
\end{align*}
and obtain similar inequalities for $I_{\text{nonlinear}}^2, I_{\text{nonlinear}}^3$. We apply the techniques used to estimate $I_{\text{linear}}^1, I_{\text{linear}}^2$ and $I_{\text{linear}}^3$ to obtain
\begin{equation*}
 I_{\text{nonlinear}} \le C \Big(\sup_{s \in [0,\delta]} \big\lVert \psi \big\rVert_{{\dot{B}}_{p,\infty}^{-1+\frac{2}{p}}} \Big) \Big(\sup_{s \in [0,\delta]} \big\lVert B(u \theta) \big\rVert_{L^{2}} \Big) \Big(\sup_{\Vert \varphi \Vert_{{\dot{B}}_{p',1}^{1-\frac{2}{p}}} = 1}  \int_{0}^{t} \big\lVert e^{(t-s)\Delta} \nabla \varphi \big\rVert_{{\dot{B}}_{p',1}^{2-\frac{2}{p}}} ds \Big).\\
\end{equation*}
Using maximum regularity \eqref{maximum regularity} for $\varphi$, we get
\begin{align*}
 I_{\text{nonlinear}} \le C \Big(\sup_{s \in [0,\delta]} \big\lVert \psi \big\rVert_{{\dot{B}}_{p,\infty}^{-1+\frac{2}{p}}} \Big) \Big(\sup_{s \in [0,\delta]} \big\lVert B(u \theta) \big\rVert_{L^{2}} \Big).
\end{align*}

Fix $\varepsilon > 0$. Then, 
applying Proposition~\ref{linear} to \eqref{liner I_1,I_3} and \eqref{linear I_2}, we get 
\begin{align}\notag
 I_{\text{linear}} \le C \varepsilon \sup_{s \in [0,\delta]} \lVert \psi \rVert_{{\dot{B}}_{p,\infty}^{-1 + \frac{2}{p}}} \ \text{for some} \ 0<\delta \ll 1.
\end{align} 
By Proposition~\ref{nonlinear} for $B(u\theta)$, 
\begin{equation}\notag
 I_{\text{nonlinear}} \le C \varepsilon \sup_{t \in [0,\delta]} \lVert \psi \rVert_{{\dot{B}}_{p,\infty}^{-1 + \frac{2}{p}}} \ \text{for some} \ 0<\delta \ll 1.
\end{equation}

Hence, 
\begin{equation}\notag
 \lVert \psi \rVert_{{\dot{B}}_{p,\infty}^{-1 + \frac{2}{p}}} \le C \varepsilon \sup_{t \in [0,\delta]} \lVert \psi \rVert_{{\dot{B}}_{p,\infty}^{-1 + \frac{2}{p}}},
\end{equation}
and as a consequence, 
\begin{equation*}
 \lVert \psi \rVert_{{\dot{B}}_{p,\infty}^{-1 + \frac{2}{p}}} \le 0 \ \text{for} \ t \in [0,\delta]
\end{equation*}
allowing us to obtain 
\begin{equation}\notag
 \theta = \tilde{\theta} \ \text{in} \ L^2 \ \text{for} \ t \in [0,\delta]
\end{equation}
by using Lemma~\ref{Besov embedding} and Lemma~\ref{shiage}.\\
\hspace{1cm}

\noindent
Step 2. In this step, we prove that $\big\lVert \psi \big\rVert_{L^2} = 0$ in the entire interval $[0,T]$, by the contradiction argument. Define
\begin{equation}\notag
 \tau^* =  \sup \big\{ \tau \in [0,T) \ \big| \ \lVert \theta(t, \cdot) - \tilde{\theta}(t, \cdot) \rVert_{L^2} = 0 \ \text{for} \ t \in [0, \tau] \big\}.
\end{equation}
If $\tau^* = T$, then the proof is completed. Assume that $\tau^* < T$. From the continuity in time of $\theta$ and $\tilde{\theta}$, it follows that $\theta(\tau^*) = \tilde{\theta}(\tau^*)$. By the same method as Step 1, the uniqueness holds in $[\tau^*, \tau^*+\delta']$ for some $\delta'>0$. This contradicts the definition of $\tau^*$.
 
 \vspace{10mm}
 \noindent
{\bf Data availability statement}. This manuscript has no associated data.

\noindent
{\bf Conflict of Interest. }
The author declares that he has no conflict of interest.

\end{document}